\documentclass[11pt]{amsart}

\usepackage[T1]{fontenc}
\usepackage[utf8]{inputenc}
\usepackage{lmodern}
\usepackage{microtype}
\usepackage[hmargin=3cm,vmargin=2.5cm]{geometry}
\usepackage{amsmath,amssymb,amsthm,mathtools}
\usepackage{xcolor}
\usepackage[colorlinks=true,linkcolor=green!50!black,citecolor=blue,urlcolor=black]{hyperref}
\usepackage{tikz}
\usepackage{float}

\theoremstyle{plain}
\newtheorem{theorem}{Theorem}[section]
\newtheorem{lemma}[theorem]{Lemma}
\newtheorem{proposition}[theorem]{Proposition}
\newtheorem{corollary}[theorem]{Corollary}

\theoremstyle{definition}
\newtheorem{example}[theorem]{Example}

\theoremstyle{remark}
\newtheorem{remark}[theorem]{Remark}

\DeclareMathOperator{\Id}{Id}
\DeclareMathOperator{\diag}{diag}
\newcommand{\RR}{\mathbb R}
\newcommand{\cC}{\mathcal C}

\title[A global rigidity criterion for Coxeter $3$-polytopes]
{A global rigidity criterion for Coxeter $3$-polytopes}
\author{Sunghwan Ko}
\address{Department of Mathematical Sciences, Seoul National University, Seoul 08826, South Korea}
\email{sym2346@snu.ac.kr}

\begin{document}

\begin{abstract}
	We give a combinatorial condition ensuring global projective rigidity of simple hyperbolic Coxeter $3$-polytopes. As an application, we prove that every right-angled L\"obell polytope is globally rigid.
\end{abstract}

\maketitle

\section{Introduction}

Let $\mathcal G$ be a labeled $3$-polytope. A projective Coxeter realization is a convex projective polytope with the prescribed combinatorics and facet reflections realizing the edge labels. We denote by $\cC(\mathcal G)$ the deformation space of these realizations modulo projective equivalence. We call $\mathcal G$ \emph{globally rigid} if $\cC(\mathcal G)$ consists of a single point.

\subsection*{Main result}

Our main result gives a combinatorial condition ensuring global rigidity for a simple labeled $3$-polytope that admits a hyperbolic realization.

Let $\mathcal G$ be a simple labeled $3$-polytope, meaning that each vertex is contained in exactly three facets. Let $S$ denote its facet set, and assume that all edge labels lie in $\{2,3,\ldots\}$. For distinct facets $F,G\in S$, write $F\sim G$ if they are adjacent, $F\not\sim G$ if they are nonadjacent, and $F\perp G$ if their common edge has label $2$. Set $n=|S|-4$ and fix a facet ordering
\[
	F_0,F_1,F_2,F_3,T_1,\ldots,T_n.
\]
For $1\le j\le n$, set
\[
	S_{j-1}=\{F_0,F_1,F_2,F_3,T_1,\ldots,T_{j-1}\}.
\]
Assume that
\begin{equation}\tag{H1}\label{hyp:initial}
	F_0\not\sim F_i\quad(1\le i\le3),
	\qquad F_i\perp F_k\quad(1\le i<k\le3),
\end{equation}
and that, for each $1\le j\le n$, there exist distinct facets $P_j,Q_j,R_j\in S_{j-1}$ satisfying
\begin{equation}\tag{H2}\label{hyp:extension}
	T_j\perp P_j,\qquad T_j\perp Q_j,\qquad T_j\perp R_j.
\end{equation}

\begin{theorem}\label{thm:criterion}
	Under assumptions~\eqref{hyp:initial} and~\eqref{hyp:extension}, every projective Coxeter polytope realizing $\mathcal G$ is hyperbolic, and $\cC(\mathcal G)$ contains at most one point. In particular, if $\mathcal G$ admits a hyperbolic realization, then it is globally rigid.
\end{theorem}

The above theorem does not guarantee the existence of a realization: the example in Appendix~\ref{app:empty-deformation-space} satisfies~\eqref{hyp:initial} and~\eqref{hyp:extension} but has an empty deformation space.

The proof of Theorem~\ref{thm:criterion} begins with the four facets in~\eqref{hyp:initial}. After a positive gauge change, their reflection data define a symmetric bilinear form of signature $(3,1)$ preserved by the four facet reflections. We then use~\eqref{hyp:extension} to prove by induction that all facet reflections preserve this form. The full Cartan matrix is then indecomposable and symmetrizable of signature $(3,1)$. We complete the proof by applying the results of Vinberg and Andreev.

\subsection*{Related results}

Global descriptions of deformation spaces are known for $2$-perfect Coxeter truncation polytopes, obtained from simplices by successive vertex truncations. Marquis studied the three-dimensional case, where $2$-perfectness means that all edge labels are finite \cite{Mar10}. Choi, Lee, and Marquis later treated higher dimensions \cite{CLM22}.

The study of deformation spaces of Coxeter polytopes is closely related to that of real projective structures on Coxeter orbifolds. Choi, Hodgson, and Lee studied restricted deformation spaces with fixed supporting hyperplanes \cite{CHL12}. They determined local dimensions for ideal polyhedra with all dihedral angles $\pi/3$ (excluding the tetrahedron) and for certain compact cubes and dodecahedra. For compact weakly orderable hyperbolic Coxeter $3$-orbifolds, Choi and Lee proved that a neighborhood of the hyperbolic structure in the projective deformation space is a cell of dimension $e_+-3$ \cite{CL15}. Here $e_+$ denotes the number of edges with labels greater than $2$.

At the infinitesimal level, Hwang and Lee \cite{HLprep} study compact hyperbolic Coxeter orbifolds that are not weakly orderable. They obtain upper bounds for the dimension of the first group cohomology associated with the holonomy. This cohomology describes infinitesimal deformations modulo infinitesimal conjugation. Using these bounds, they establish local projective rigidity for suitable classes of these orbifolds.

\subsection*{Applications}

Our main application is to the L\"obell polytopes $L_m$, $m\ge5$. Each $L_m$ has two $m$-gonal facets and two belts of $m$ pentagonal facets. For example, $L_5$ is the dodecahedron. These polytopes are not truncation polytopes, since they have no triangular facets. The global descriptions of deformation spaces for truncation polytopes discussed above do not apply to this family.

Andreev's theorem \cite{And1970a,And1970b} gives a compact hyperbolic realization of the right-angled labeling of each $L_m$. This labeling admits an ordering satisfying~\eqref{hyp:initial} and~\eqref{hyp:extension}, as shown in Proposition~\ref{prop:lobell-labelings}. Theorem~\ref{thm:criterion} therefore gives the following corollary.

\begin{corollary}\label{cor:main}
	For every $m\ge5$, the deformation space of the right-angled L\"obell polytope $L_m$ consists of a single point.
\end{corollary}

More generally, Proposition~\ref{prop:lobell-labelings} characterizes the finite edge labelings of $L_m$ that satisfy~\eqref{hyp:initial} and~\eqref{hyp:extension}. Each of these labelings admits a hyperbolic realization by~\cite{Atk2011} and is therefore globally rigid; see Corollary~\ref{cor:lobell-labelings}. This gives examples with non-right edges and arbitrarily many facets, including finite, ideal, or hyperideal vertices; see Example~\ref{ex:three-edges}.

\subsection*{Organization}
Section~\ref{sec:preliminaries} reviews the reflection data and the results of Vinberg and Andreev. Section~\ref{sec:proof} proves Theorem~\ref{thm:criterion}. Section~\ref{sec:lobell} characterizes the L\"obell labelings and proves the rigidity corollaries. Appendix~\ref{app:empty-deformation-space} gives an explicit example with an empty deformation space. Appendix~\ref{app:non-lobell-example} gives a non-L\"obell example satisfying the criterion.

\subsection*{Acknowledgements}
I am grateful to my advisor, Gye-Seon Lee, for many helpful discussions on this project. This work was supported by the National Research Foundation of Korea(NRF) grant funded by the Korea government(MSIT) (RS-2023-00252171).

\subsection*{AI disclosure}
GPT Sol~5.6 suggested the idea underlying Lemma~\ref{lem:extension} during discussions on the rigidity of a right-angled dodecahedron. The author developed the argument and formulated conditions~\eqref{hyp:initial} and~\eqref{hyp:extension} through further discussions with GPT-6 Astra. AI tools also assisted with editing the text and preparing the TikZ figures. The author takes full responsibility for the paper.

\section{Preliminaries}
\label{sec:preliminaries}

Let $V=\RR^4$ and $\mathbb P(V)=(V\setminus\{0\})/\RR^*$. For $U\subset V$, let $\mathbb P(U)$ denote the image of $U\setminus\{0\}$. A subset of $\mathbb P(V)$ is \emph{convex} if it is the projectivization of a convex cone. A convex subset is \emph{properly convex} if its closure is contained in an affine chart.
A \emph{projective $3$-polytope} is a properly convex subset $P\subset\mathbb P(V)$ with nonempty interior that can be written as
\[
	P=\mathbb P(C_P),\qquad C_P=\bigcap_{s\in S}\{x\in V:\alpha_s(x)\le0\},
\]
where $0\neq\alpha_s\in V^*$ and the defining inequalities are irredundant. The facets are indexed by $S$, with $F_s$ supported by $\mathbb P(\ker\alpha_s)$.

A projective reflection with support $\mathbb P(\ker\alpha_s)$ and pole $[b_s]\in \mathbb P(V)$ has a linear lift
\[
	\sigma_s=\Id-\alpha_s\otimes b_s,\qquad
	\alpha_s(b_s)=2.
\]
The vector $b_s$ is called a \emph{pole vector}. The \emph{positive gauge action} $(\alpha_s,b_s)\mapsto(\lambda_s\alpha_s,\lambda_s^{-1}b_s)$ for $\lambda_s>0$ preserves both $\sigma_s$ and the defining half-space $\{x\in V:\alpha_s(x)\le0\}$.

Let $\mathcal G$ be a labeled $3$-polytope with facet set $S$ and edge labels $m_e\in\{2,3,\ldots\}$. For distinct $s,t\in S$, let $m_{st}=m_e$ if $s$ and $t$ meet along an edge $e$, and set $m_{st}=\infty$ otherwise. Let $\mathcal F(Q)$ denote the face poset of a polytope $Q$.
A \emph{Coxeter polytope realizing $\mathcal G$} is a pair $(P,\phi)$, where $P$ is a projective $3$-polytope equipped with a reflection in each facet and $\phi:\mathcal F(\mathcal G)\to\mathcal F(P)$ is a label-preserving poset isomorphism. Writing $F_s=\phi(s)$, we require that, for all distinct $s,t\in S$, the associated reflection data satisfy $\alpha_s(b_t)\le0$ together with
\begin{equation}\label{eq:cartan equations}
	\begin{aligned}
		\alpha_s(b_t)=\alpha_t(b_s)=0
		&\quad\text{if }m_{st}=2,\\
		\alpha_s(b_t)\alpha_t(b_s)=4\cos^2\frac{\pi}{m_{st}}&\quad \text{if }3\le m_{st}<\infty,\\
		\alpha_s(b_t)\alpha_t(b_s)\ge4
		&\quad\text{if }m_{st}=\infty.
	\end{aligned}
\end{equation}

Two Coxeter polytopes $(P,\phi)$ and $(P',\phi')$ realizing $\mathcal G$ are \emph{isomorphic} if a projective automorphism $\psi$ of $\mathbb P(V)$ maps $P$ onto $P'$, conjugates the corresponding facet reflections, and satisfies $\widehat\psi\circ\phi=\phi'$, where $\widehat\psi$ denotes the induced map of face posets. We omit the marking $\phi$ when no confusion can arise. The deformation space $\cC(\mathcal G)$ is the set of isomorphism classes of Coxeter polytopes realizing $\mathcal G$.

The \emph{Cartan matrix} of $P$ is $A_P=(a_{st})_{s,t\in S}$ with $a_{st}=\alpha_s(b_t)$. The positive gauge action induces
\[
	A_P\longmapsto\Lambda A_P\Lambda^{-1}, \qquad \Lambda=\diag(\lambda_s)_{s\in S}.
\]
The matrix $A_P$ is \emph{symmetrizable} if $\Lambda A_P\Lambda^{-1}$ is symmetric for some positive diagonal matrix $\Lambda$. Its signature is that of this symmetric representative with zero eigenvalues omitted. It is \emph{indecomposable} if the graph on $S$ joining distinct $s,t$ whenever $a_{st}\ne0$ is connected.

A nondegenerate symmetric bilinear form $B$ on $V$ of signature $(3,1)$ defines the Klein model
\[
	\mathbb H_B^3=\{[x]\in\mathbb P(V):B(x,x)<0\}.
\]
A Coxeter polytope $P$ is \emph{hyperbolic} if its facet reflections act as hyperbolic reflections on some $\mathbb H_B^3$ that meets the interior of $P$. Its \emph{hyperbolic chamber} $Q_P=P\cap\mathbb H_B^3$ has the prescribed dihedral angles by~\eqref{eq:cartan equations}. A vertex of $P$ is \emph{finite}, \emph{ideal}, or \emph{hyperideal} according as it lies in $\mathbb H_B^3$, in $\partial\mathbb H_B^3$, or outside $\overline{\mathbb H_B^3}$.

The following result of Vinberg relates symmetrizable Cartan matrices to hyperbolic Coxeter polytopes.

\begin{theorem}[Vinberg {\cite{Vin71}}]
	\label{thm:vinberg-cartan}
	Let $P$ be a Coxeter $3$-polytope. If its Cartan matrix $A_P$ is indecomposable and symmetrizable of signature $(3,1)$, then $P$ is hyperbolic.
\end{theorem}

A $3$-polytope is \emph{simple} if each vertex is contained in exactly three facets. A \emph{prismatic $k$-circuit} is a simple $k$-cycle in the dual $1$-skeleton whose crossed primal edges have pairwise disjoint sets of endpoints.

We use the following consequences of Andreev's theorems and Atkinson's extension to generalized hyperbolic polyhedra. See also Bao and Bonahon \cite{BB2002} for a characterization of hyperideal polyhedra.

\begin{theorem}[Andreev--Atkinson {\cite{And1970a,And1970b,Atk2011}}]
	\label{thm:hyperbolic-existence}
	Let $\mathcal G$ be a simple labeled $3$-polytope with finite edge labels $m_e$, and set $\theta_e=\pi/m_e$. If $\mathcal G$ admits a hyperbolic realization, then the following conditions hold:
	\begin{enumerate}
		\item For every prismatic $3$-circuit crossing edges $e_1,e_2,e_3$,
		\[
			\theta_{e_1}+\theta_{e_2}+\theta_{e_3}<\pi.
		\]
		\item For every prismatic $4$-circuit crossing edges $e_1,e_2,e_3,e_4$,
		\[
			\theta_{e_1}+\theta_{e_2}+\theta_{e_3}+\theta_{e_4}<2\pi.
		\]
	\end{enumerate}
	If $\mathcal G$ has more than four vertices and no prismatic $3$- or $4$-circuits, then it admits a hyperbolic realization. In any hyperbolic realization, a vertex with incident labels $p,q,r$ is finite, ideal, or hyperideal according as $1/p+1/q+1/r$ is greater than, equal to, or less than $1$.
\end{theorem}

\begin{theorem}[Andreev {\cite{And1970a,And1970b}}]
	\label{thm:hyperbolic-uniqueness}
	Let $\mathcal G$ be a simple labeled $3$-polytope with finite edge labels. Any two hyperbolic Coxeter polytopes realizing $\mathcal G$ are isomorphic.
\end{theorem}

For realizations with hyperideal vertices, the uniqueness statement follows by truncating those vertices along their polar planes and applying Andreev's finite-volume theorem. The new dihedral angles are all $\pi/2$; see~\cite[Section~5]{Atk2011}. The original facet planes are then determined by the truncated polyhedron.

\section{Proof of Theorem~\ref{thm:criterion}}\label{sec:proof}

Throughout this section, let $P$ be a Coxeter $3$-polytope with reflection data $(\alpha_F,b_F)_{F\in S}$ on $V=\RR^4$ and Cartan matrix $A=A_P=(a_{FG})_{F,G\in S}$.

\begin{lemma}\label{lem:three-covectors}
	The defining covectors of any three distinct facets of a properly convex polytope with nonempty interior are linearly independent.
\end{lemma}

\begin{proof}
	Let $\alpha_1,\alpha_2,\alpha_3$ define three distinct facets. Since all defining covectors are negative on the interior of the defining cone, any nontrivial dependence has coefficients of both signs. After relabeling, it gives $\alpha_1=u\alpha_2+v\alpha_3$ with $u,v\ge0$ and $u+v>0$. Choose $x$ in the defining cone projecting to a relative interior point of the first facet. Then $\alpha_1(x)=0$ while $\alpha_2(x),\alpha_3(x)<0$, a contradiction.
\end{proof}

\begin{lemma}\label{lem:initial-four}
	Suppose that distinct facets $F_0,F_1,F_2,F_3$ satisfy~\eqref{hyp:initial}. After a positive gauge change, the corresponding pole vectors and covectors form bases of $V$ and $V^*$, respectively. There is a unique symmetric bilinear form $B$ of signature $(3,1)$ satisfying
	\begin{equation}\label{eq:initial-duality}
		B(b_{F_i},\cdot)=\alpha_{F_i}
		\qquad(0\le i\le3).
	\end{equation}
\end{lemma}

\begin{proof}
	Apply the positive gauge with $\lambda_{F_i}=\sqrt{a_{F_0F_i}/a_{F_iF_0}}$ for $1\le i\le3$ and all other factors equal to $1$. The Cartan conditions and~\eqref{hyp:initial} give the initial principal submatrix
	\[
		A_0=
		\begin{pmatrix}
			2&-c_1&-c_2&-c_3\\
			-c_1&2&0&0\\
			-c_2&0&2&0\\
			-c_3&0&0&2
		\end{pmatrix},
		\qquad
		c_i=\sqrt{a_{F_0F_i}a_{F_iF_0}}\ge2.
	\]
	The Schur complement of $2I_3$ is $2-\frac12\sum_{i=1}^3c_i^2<0$, so $A_0$ is nonsingular with signature $(3,1)$. The initial pole vectors and covectors therefore form bases of $V$ and $V^*$, respectively. The form $B$ with Gram matrix $A_0$ in the basis of initial pole vectors satisfies~\eqref{eq:initial-duality} and is uniquely determined by it.
\end{proof}

\begin{lemma}\label{lem:extension}
	Let $\mathcal S\subseteq S$ satisfy
	\[
		\operatorname{span}\{\alpha_F:F\in\mathcal S\}=V^*.
	\]
	Suppose that a nondegenerate symmetric bilinear form $B$ on $V$ and positive numbers $(d_F)_{F\in\mathcal S}$ satisfy
	\[
		B(b_F,\cdot)=d_F\alpha_F
		\qquad(F\in\mathcal S).
	\]
	If $T\in S\setminus\mathcal S$ is orthogonal to three distinct facets in $\mathcal S$, then there is a unique $d_T>0$ such that
	\[
		B(b_T,\cdot)=d_T\alpha_T.
	\]
\end{lemma}

\begin{proof}
	Let $K_1,K_2,K_3\in\mathcal S$ be distinct facets orthogonal to $T$. By the assumed duality and Lemma~\ref{lem:three-covectors}, their pole vectors span a three-dimensional subspace $W$. Orthogonality implies that both $\alpha_T$ and $B(b_T,\cdot)$ vanish on $W$. These nonzero covectors therefore lie in the one-dimensional annihilator of $W$, so $B(b_T,\cdot)=d_T\alpha_T$ for a unique $d_T\ne0$.

	The spanning hypothesis gives $E\in\mathcal S$ with $a_{ET}\ne0$. The Cartan conditions~\eqref{eq:cartan equations} imply $a_{ET},a_{TE}<0$, and symmetry yields
	\[
		d_T=d_E\frac{a_{ET}}{a_{TE}}>0.
	\]
\end{proof}

\begin{proof}[Proof of Theorem~\ref{thm:criterion}]
	Let $P$ realize $\mathcal G$. After a positive gauge change, Lemma~\ref{lem:initial-four} gives a form $B$ of signature $(3,1)$ satisfying~\eqref{eq:initial-duality}. Since the initial covectors span $V^*$, successive applications of Lemma~\ref{lem:extension} using~\eqref{hyp:extension} give
	\begin{equation}\label{eq:global-duality}
		B(b_F,\cdot)=d_F\alpha_F,
		\qquad d_F>0,
		\qquad F\in S.
	\end{equation}
	Set $\Lambda=\diag(\sqrt{d_F})_{F\in S}$. The matrix $\Lambda A\Lambda^{-1}$ is the Gram matrix of the spanning family $(b_F/\sqrt{d_F})_{F\in S}$ with respect to $B$. Thus $A$ is symmetrizable of signature $(3,1)$.

	The initial four indices lie in one connected component of the graph of $A$ by~\eqref{hyp:initial}. Every remaining column has a nonzero entry in one of the initial four rows since the initial covectors span $V^*$ and each pole vector is nonzero. Thus $A$ is indecomposable.

	By Theorems~\ref{thm:vinberg-cartan} and~\ref{thm:hyperbolic-uniqueness}, $P$ is hyperbolic and $\cC(\mathcal G)$ contains at most one point.
\end{proof}

\begin{remark}\label{rem:unused-edges}
	Fix the facet ordering and the triples in~\eqref{hyp:extension}, and let $N=|S|$. The three edges among $F_1,F_2,F_3$ and the three selected for each $T_j$ give $3+3(N-4)=3N-9$ distinct edges labeled $2$. No edge is counted twice, since an edge in~\eqref{hyp:extension} is selected only when its later incident facet is added. Since $\mathcal G$ is simple, Euler's formula gives $3N-6$ edges in total. Exactly three edges are therefore unselected, so at most three edges have labels greater than $2$.

	Replacing the labels on the unselected edges by arbitrary finite integers at least $2$ preserves~\eqref{hyp:initial} and~\eqref{hyp:extension}. Every resulting labeling that admits a realization is therefore globally rigid.
\end{remark}

\section{Applications to labeled L\"obell polytopes}\label{sec:lobell}

Fix $m\ge5$, and interpret all indices modulo $m$. The facets of the L\"obell polytope $L_m$ are
\[
	U,D,X_0,\ldots,X_{m-1},Y_0,\ldots,Y_{m-1},
\]
where $U,D$ are $m$-gons and all remaining facets are pentagons. The adjacency relations are
\begin{equation}\label{eq:lobell-adjacencies}
	\begin{aligned}
		U&\sim X_i, & D&\sim Y_i,\\
		X_i&\sim X_{i-1},X_{i+1},Y_{i-1},Y_i,\\
		Y_i&\sim Y_{i-1},Y_{i+1},X_i,X_{i+1}.
	\end{aligned}
\end{equation}

\begin{figure}[!htbp]
	\centering
	\begin{tikzpicture}[x=0.60cm,y=0.60cm,line width=0.45pt,line cap=round,line join=round,font=\normalsize,every node/.style={inner sep=0pt}]
		\path[use as bounding box] (-4.55,-4.8) rectangle (4.55,4.8);
		\foreach \i in {0,...,5} {
			\pgfmathsetmacro{\lobellangle}{90+60*\i}
			\coordinate (a\i) at (\lobellangle:4.5);
			\coordinate (b\i) at (\lobellangle:2.8);
			\coordinate (c\i) at ({\lobellangle-30}:2.3);
			\coordinate (d\i) at ({\lobellangle+30}:1.3);
		}
		\foreach \i in {0,...,5} {
			\pgfmathtruncatemacro{\j}{mod(\i+1,6)}
			\draw (a\i)--(a\j);
			\draw (a\i)--(b\i);
			\draw (c\i)--(b\i)--(c\j);
			\draw (c\j)--(d\i);
			\draw (d\i)--(d\j);
		}
		\foreach \ring in {a,b,c,d} {
			\foreach \i in {0,...,5} {
				\fill (\ring\i) circle[radius=0.65pt];
			}
		}
		\node at (3.25,3.85) {$U$};
		\node at (0,0) {$D$};
		\foreach \i in {0,...,5} {
			\node at ({60+60*\i}:3.23) {$X_{\i}$};
			\node at ({90+60*\i}:1.82) {$Y_{\i}$};
		}
	\end{tikzpicture}
	\caption{Facet notation for $L_6$. The unbounded face represents $U$, the central hexagon represents $D$, and indices are taken modulo $6$.}
	\label{fig:lobell-facets}
\end{figure}
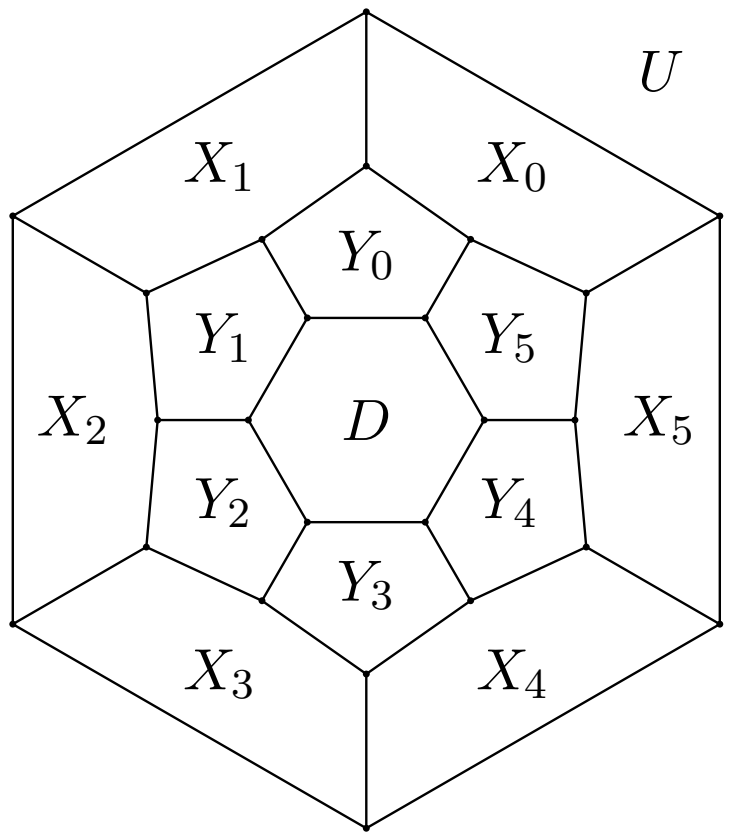

For each facet $F$, let $E(F)$ denote its set of edges. Given a labeling $\mathcal G$ of $L_m$ with finite edge labels, set $E_+=\{e:m_e>2\}$. The following proposition characterizes the labelings that satisfy the hypotheses of Theorem~\ref{thm:criterion}.

\begin{proposition}\label{prop:lobell-labelings}
	The labeled polytope $\mathcal G$ admits a facet ordering satisfying~\eqref{hyp:initial} and~\eqref{hyp:extension} if and only if $|E_+|\le3$ and there exist adjacent pentagonal facets $A,B$ such that
	\[
		E_+\subseteq E(A)\cup E(B),\qquad |E_+\cap E(A)|\le2,\qquad |E_+\cap E(B)|\le2.
	\]
\end{proposition}

\begin{proof}
    Suppose an ordering satisfies~\eqref{hyp:initial} and~\eqref{hyp:extension}. Let $E_0\supseteq E_+$ be the unselected edge set in Remark~\ref{rem:unused-edges}. Since $L_m$ is simple, $|E_0|=3$. Let $A$ and $B$ be the last and second-last facets, with $a$ and $b$ edges. Set $\varepsilon=1$ if they are adjacent and $\varepsilon=0$ otherwise. Exactly three edges joining $A$ to earlier facets are selected in~\eqref{hyp:extension}, as are exactly three joining $B$ to earlier facets. The remaining edges to earlier facets form disjoint subsets of $E_0$ of sizes $a-3$ and $b-\varepsilon-3$, respectively. Since $a,b\ge5$,
	\[
		3\ge(a-3)+(b-\varepsilon-3)\ge3.
	\]
	Thus $a=b=5$ and $\varepsilon=1$. The three edges of $E_0$ are two on $A$ and one on $B$ other than $A\cap B$. Thus $E_0\subseteq E(A)\cup E(B)$ with at most two edges on either facet. The same bounds hold for $E_+\subseteq E_0$.

	Conversely, suppose the stated conditions hold. All edges outside $E(A)\cup E(B)$ have label $2$. Using rotations, reflections and an exchange of the two belts, we may assume that $(A,B)=(X_{m-1},Y_{m-1})$ or $(A,B)=(X_0,X_1)$. When $m=5$, symmetries of the dodecahedron also cover pairs involving $U$ or $D$.

	For $(A,B)=(X_{m-1},Y_{m-1})$, take $(F_0,F_1,F_2,F_3)=(D,U,X_0,X_1)$. Then add $Y_0,\allowbreak Y_1,\allowbreak X_2,\allowbreak Y_2,\allowbreak\ldots,\allowbreak X_{m-2},\allowbreak Y_{m-2}$ using the following triples:
	\[
		\begin{array}{c|ccc}
			T_j & P_j & Q_j & R_j\\ \hline
			Y_0 & D & X_0 & X_1\\
			Y_k\ (1\le k\le m-2) & D & X_k & Y_{k-1}\\
			X_k\ (2\le k\le m-2) & U & X_{k-1} & Y_{k-1}
		\end{array}
	\]

	For $(A,B)=(X_0,X_1)$, take $(F_0,F_1,F_2,F_3)=(U,D,Y_1,Y_2)$. Then add $X_2,\allowbreak X_3,\allowbreak Y_3,\allowbreak\ldots,\allowbreak X_{m-1},\allowbreak Y_{m-1},\allowbreak Y_0$ using the following triples:
	\[
		\begin{array}{c|ccc}
			T_j & P_j & Q_j & R_j\\ \hline
			X_2 & U & Y_1 & Y_2\\
			X_k\ (3\le k\le m-1) & U & X_{k-1} & Y_{k-1}\\
			Y_k\ (3\le k\le m-1) & D & X_k & Y_{k-1}\\
			Y_0 & D & Y_1 & Y_{m-1}
		\end{array}
	\]
	Both initial choices satisfy~\eqref{hyp:initial}. Each triple consists of earlier facets and selects edges outside $E(A)\cup E(B)$, so~\eqref{hyp:extension} holds.

	It remains to add $A$ and $B$ in a suitable order. Since $|E_+|\le3$, we may choose the first facet to have at most one non-right edge other than $A\cap B$. Its four edges other than $A\cap B$ are shared with earlier facets, and at least three have label $2$. The second facet is adjacent to five earlier facets and has at most two non-right edges. Hence each of $A$ and $B$ is orthogonal to at least three earlier facets, as required by~\eqref{hyp:extension}.
\end{proof}

\begin{corollary}\label{cor:lobell-labelings}
	Let $\mathcal G$ be a labeling of $L_m$ with finite edge labels, where $m\ge5$. Suppose that at most three edges have labels greater than $2$ and that all these edges are contained in the union of the edge sets of two adjacent pentagonal facets, each of which contains at most two of these edges. Then $\cC(\mathcal G)$ consists of a single point.
\end{corollary}

\begin{proof}
	The L\"obell polytope $L_m$ is simple, has $4m>4$ vertices, and has no prismatic $3$- or $4$-circuits. Thus Theorem~\ref{thm:hyperbolic-existence} gives a hyperbolic realization of $\mathcal G$. Proposition~\ref{prop:lobell-labelings} and Theorem~\ref{thm:criterion} imply that $\cC(\mathcal G)$ consists of a single point.
\end{proof}

\begin{proof}[Proof of Corollary~\ref{cor:main}]
	Apply Corollary~\ref{cor:lobell-labelings} to the labeling of $L_m$ in which every edge has label $2$.
\end{proof}

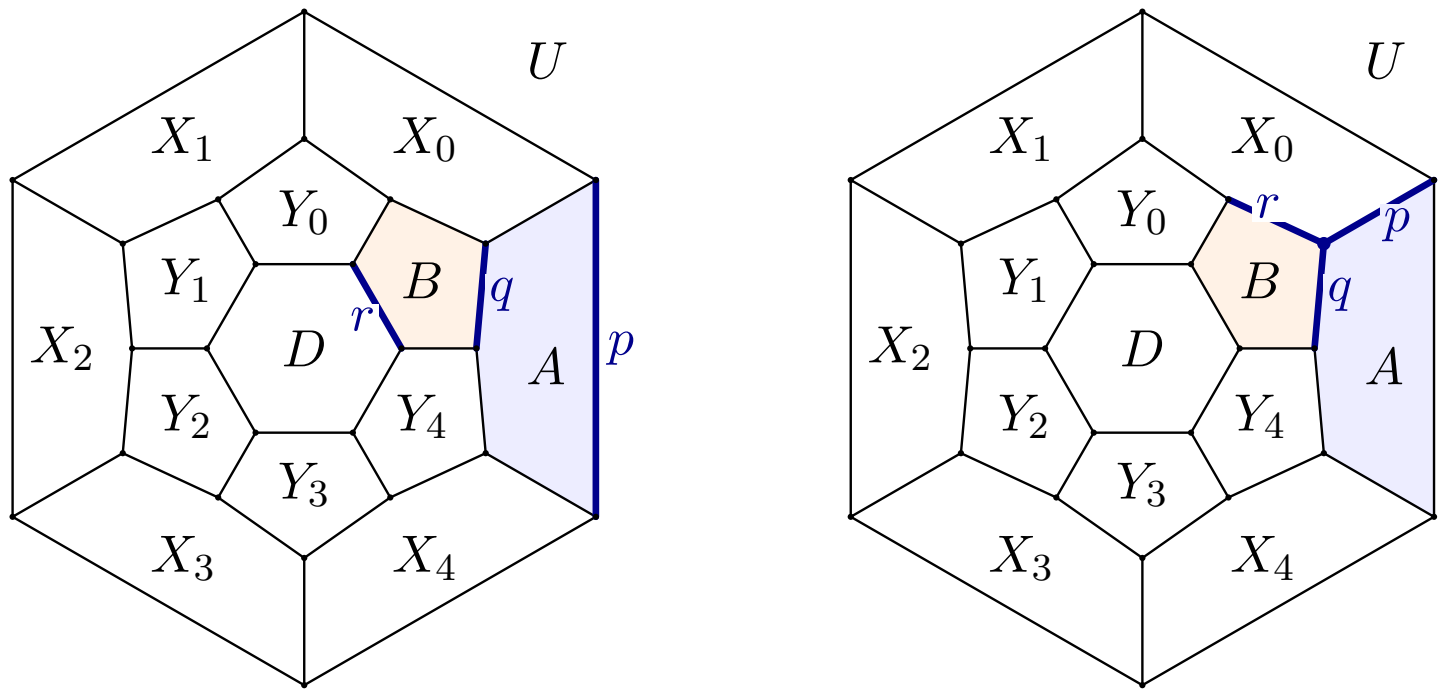
\begin{figure}[H]
	\centering
	\begin{tikzpicture}[x=0.55cm,y=0.55cm,scale=0.9,transform shape,line width=0.45pt,line cap=round,line join=round,font=\normalsize,every node/.style={inner sep=0pt}]
		\path[use as bounding box] (-4.55,-5.6) rectangle (15.75,4.8);
		\foreach \panel/\offset in {0/0,1/11.2} {
			\begin{scope}[shift={(\offset,0)}]
				\foreach \i in {0,...,5} {
					\pgfmathsetmacro{\lobellangle}{90+60*\i}
					\coordinate (a\i) at (\lobellangle:4.5);
					\coordinate (b\i) at (\lobellangle:2.8);
					\coordinate (c\i) at ({\lobellangle-30}:2.3);
					\coordinate (d\i) at ({\lobellangle+30}:1.3);
				}
				\fill[blue!7] (a4)--(a5)--(b5)--(c5)--(b4)--cycle;
				\fill[orange!10] (c5)--(b5)--(c0)--(d5)--(d4)--cycle;
				\foreach \i in {0,...,5} {
					\pgfmathtruncatemacro{\j}{mod(\i+1,6)}
					\draw (a\i)--(a\j);
					\draw (a\i)--(b\i);
					\draw (c\i)--(b\i)--(c\j);
					\draw (c\j)--(d\i);
					\draw (d\i)--(d\j);
				}
				\ifnum\panel=0
					\draw[blue!55!black,line width=1.25pt] (a4)--(a5);
					\draw[blue!55!black,line width=1.25pt] (b5)--(c5);
					\draw[blue!55!black,line width=1.25pt] (d4)--(d5);
				\else
					\draw[blue!55!black,line width=1.25pt] (a5)--(b5);
					\draw[blue!55!black,line width=1.25pt] (b5)--(c5);
					\draw[blue!55!black,line width=1.25pt] (b5)--(c0);
				\fi
				\foreach \ring in {a,b,c,d} {
					\foreach \i in {0,...,5} {
						\fill (\ring\i) circle[radius=0.65pt];
					}
				}
				\node at (3.25,3.85) {$U$};
				\node at (0,0) {$D$};
				\foreach \i in {0,...,4} {
					\node at ({60+60*\i}:3.23) {$X_{\i}$};
					\node at ({90+60*\i}:1.82) {$Y_{\i}$};
				}
				\node at (3.23,-0.24) {$A$};
				\node at (30:1.82) {$B$};
				\node[text=blue!55!black,fill=blue!7,inner sep=0.8pt] at (2.64,0.72) {$q$};
				\ifnum\panel=0
					\node[text=blue!55!black,anchor=west] at (4.04,0) {$p$};
					\node[text=blue!55!black,fill=white,inner sep=0.8pt] at (0.78,0.39) {$r$};
					\node at (0,-5.2) {(a) Three disjoint edges};
				\else
					\fill[blue!55!black] (b5) circle[radius=1.4pt];
					\node[text=blue!55!black,fill=blue!7,inner sep=0.8pt] at (3.4,1.65) {$p$};
					\node[text=blue!55!black,fill=white,inner sep=0.8pt] at (1.68,1.9) {$r$};
					\node at (0,-5.2) {(b) Three edges at one vertex};
				\fi
			\end{scope}
		}
	\end{tikzpicture}
	\caption{Two labelings of $L_6$ from Example~\ref{ex:three-edges}, with $A=X_5$ and $B=Y_5$. The highlighted edges carry arbitrary finite labels $p,q,r\ge2$; every other edge has label $2$.}
	\label{fig:lobell-pqr}
\end{figure}

\begin{example}\label{ex:three-edges}
	Let $m\ge5$, and set $A=X_{m-1}$ and $B=Y_{m-1}$ in $L_m$. Label $U\cap A$, $A\cap B$, and $D\cap B$ by arbitrary finite integers $p,q,r\ge2$, respectively, and all other edges by $2$; see Figure~\ref{fig:lobell-pqr}(a). Corollary~\ref{cor:lobell-labelings} gives global rigidity, and every vertex is finite since its incident labels are $(2,2,k)$ with $k<\infty$.

	Alternatively, label the three edges at one vertex by $p,q,r\ge2$ and all other edges by $2$; see Figure~\ref{fig:lobell-pqr}(b). Corollary~\ref{cor:lobell-labelings} again gives global rigidity. The chosen vertex is finite, ideal, or hyperideal according as $1/p+1/q+1/r$ is greater than, equal to, or less than $1$; all other vertices are finite.
\end{example}

\clearpage
\appendix
\section{An empty deformation space}\label{app:empty-deformation-space}
\setcounter{figure}{0}
\renewcommand{\thefigure}{\thesection.\arabic{figure}}
\renewcommand{\theHfigure}{\thesection.\arabic{figure}}

The hypotheses of Theorem~\ref{thm:criterion} do not ensure the existence of a realization. We verify this explicitly for a vertex truncation of a dodecahedron.

\begin{figure}[H]
	\centering
	\begin{tikzpicture}[x=0.55cm,y=0.55cm,scale=0.9,transform shape,line width=0.45pt,line cap=round,line join=round,font=\normalsize,every node/.style={inner sep=0pt}]
		\path[use as bounding box] (-4.65,-4.65) rectangle (15.65,5.45);
		\foreach \panel/\offset in {0/0,1/11} {
			\begin{scope}[shift={(\offset,0)}]
				\foreach \i in {0,...,4} {
					\pgfmathsetmacro{\dodecaangle}{90+72*\i}
					\coordinate (a\i) at (\dodecaangle:4.5);
					\coordinate (b\i) at (\dodecaangle:2.8);
					\coordinate (c\i) at ({\dodecaangle-36}:2.3);
					\coordinate (d\i) at ({\dodecaangle+36}:1.3);
				}
				\path (a0)--(a1) coordinate[pos=0.12] (tL);
				\path (a0)--(a4) coordinate[pos=0.12] (tR);
				\path (a0)--(b0) coordinate[pos=0.75] (tB);
				\ifnum\panel=1
					\fill[orange!20] (tL)--(tR)--(tB)--cycle;
				\fi
				\foreach \i in {0,...,4} {
					\pgfmathtruncatemacro{\j}{mod(\i+1,5)}
					\ifnum\panel=0
						\draw (a\i)--(a\j);
						\draw (a\i)--(b\i);
					\else
						\ifnum\i=0
							\draw (tL)--(a1);
							\draw (tB)--(b0);
						\else
							\ifnum\i=4
								\draw (a4)--(tR);
							\else
								\draw (a\i)--(a\j);
							\fi
							\draw (a\i)--(b\i);
						\fi
					\fi
					\draw (c\i)--(b\i)--(c\j);
					\draw (c\j)--(d\i);
					\draw (d\i)--(d\j);
				}
				\ifnum\panel=0
					\draw[blue!55!black,line width=1.05pt] (a1)--(a0)--(a4);
					\draw[blue!55!black,line width=1.05pt] (a0)--(b0);
					\fill[blue!55!black] (a0) circle[radius=1.1pt];
					\node[above=3pt,text=blue!55!black] at (a0) {$v$};
				\else
					\draw[blue!55!black,line width=1.05pt] (a1)--(tL);
					\draw[blue!55!black,line width=1.05pt] (a4)--(tR);
					\draw[blue!55!black,line width=1.05pt] (b0)--(tB);
					\draw (tL)--(tR)--(tB)--cycle;
					\foreach \vertex in {tL,tR,tB} {
						\fill (\vertex) circle[radius=0.65pt];
					}
					\draw[red!65!black,dash pattern=on 2pt off 1.5pt,line width=0.65pt]
						(0,4.08) ellipse [x radius=1.65,y radius=1.18];
					\node[font=\small] at (0,3.88) {$T$};
				\fi
				\foreach \ring in {b,c,d} {
					\foreach \i in {0,...,4} {
						\fill (\ring\i) circle[radius=0.65pt];
					}
				}
				\foreach \i in {1,...,4} {
					\fill (a\i) circle[radius=0.65pt];
				}
				\node at (3.3,4.2) {$U$};
				\node at (0,0) {$D$};
				\foreach \i in {0,...,4} {
					\node at ({54+72*\i}:3.10) {$X_{\i}$};
					\node at ({90+72*\i}:1.82) {$Y_{\i}$};
				}
				\ifnum\panel=0
					\node at (0,-4.35) {(a) Dodecahedron};
				\else
					\node at (0,-4.35) {(b) Vertex truncation};
				\fi
			\end{scope}
		}
	\end{tikzpicture}
	\caption{A vertex truncation of the dodecahedron at $v=U\cap X_0\cap X_1$. The shaded facet is $T$. In~(b), the dashed curve crosses the three highlighted old edges and represents a prismatic $3$-circuit. The unbounded face is $U$ in both panels, and every edge has label $2$.}
	\label{fig:truncated-dodecahedron}
\end{figure}
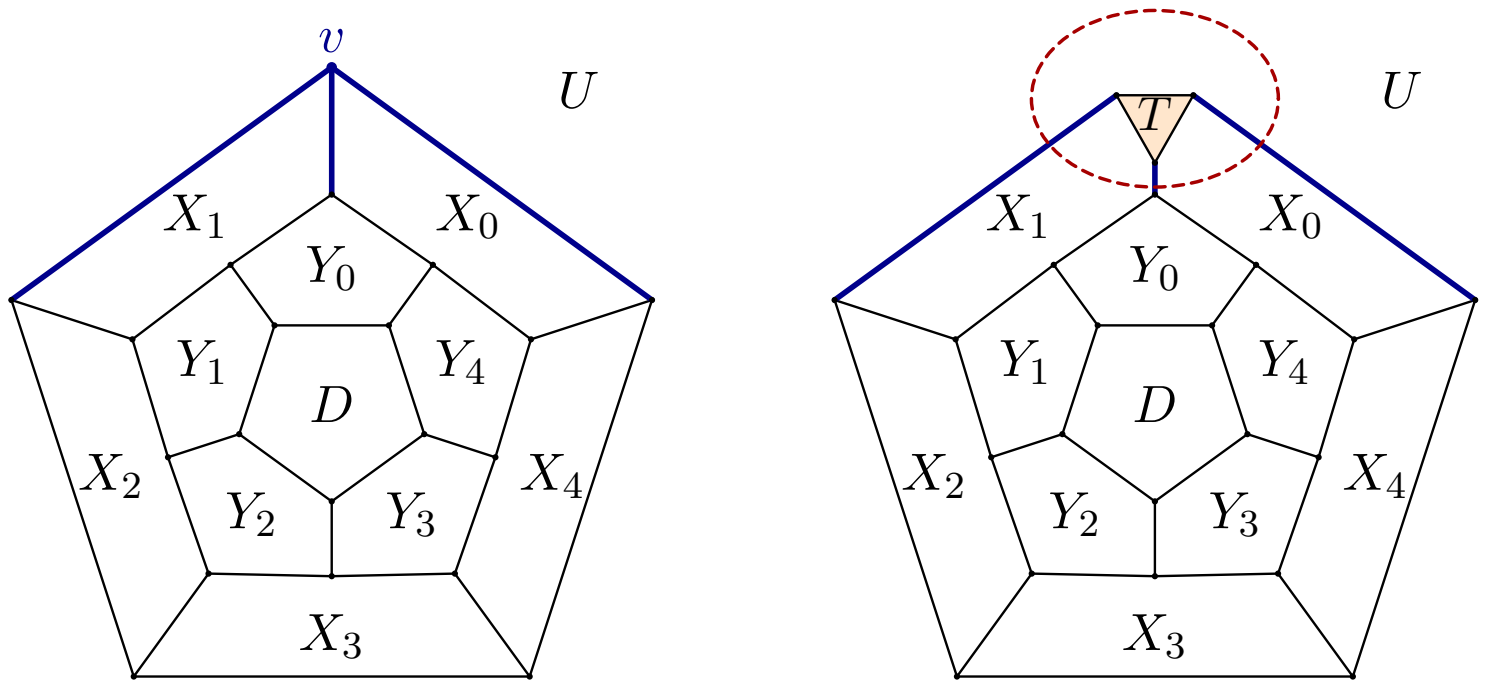

Label the facets of a dodecahedron by $U,D,X_0,\ldots,X_4,Y_0,\ldots,Y_4$ as in Figure~\ref{fig:truncated-dodecahedron}(a). Truncate the vertex $v=U\cap X_0\cap X_1$, denote the new triangular facet by $T$, and retain the names of the old facets. Let $\mathcal G'$ be the resulting polytope with every edge labeled $2$; see Figure~\ref{fig:truncated-dodecahedron}(b).

Vertex truncation preserves all adjacencies among the old facets. Thus the initial choice
\[
	(F_0,F_1,F_2,F_3)=(D,U,X_0,X_1)
\]
satisfies~\eqref{hyp:initial}: $D$ is nonadjacent to $U,X_0,X_1$, and these three facets are pairwise orthogonal. Order the remaining facets as
\[
	Y_0,Y_1,X_2,Y_2,X_3,Y_3,X_4,Y_4,T.
\]
The following choices of earlier orthogonal facets verify~\eqref{hyp:extension}:
\[
	\begin{array}{c|ccc}
		T_j & P_j & Q_j & R_j\\ \hline
		Y_0 & D & X_0 & X_1\\
		Y_k\ (1\le k\le4) & D & X_k & Y_{k-1}\\
		X_k\ (2\le k\le4) & U & X_{k-1} & Y_{k-1}\\
		T & U & X_0 & X_1
	\end{array}
\]

However, the three old facets $U,X_0,X_1$ form a prismatic $3$-circuit after truncation. The dashed curve in Figure~\ref{fig:truncated-dodecahedron}(b) crosses the shortened old edges $U\cap X_0$, $U\cap X_1$, and $X_0\cap X_1$, whose endpoints are pairwise disjoint. The corresponding dihedral angles sum to $3\pi/2$, violating the prismatic $3$-circuit inequality in Theorem~\ref{thm:hyperbolic-existence}. Thus $\mathcal G'$ has no hyperbolic realization, and Theorem~\ref{thm:criterion} implies that $\cC(\mathcal G')=\varnothing$.

\clearpage
\section{An example of a non-L\"obell polytope}\label{app:non-lobell-example}
\setcounter{figure}{0}
	Conditions~\eqref{hyp:initial} and~\eqref{hyp:extension} hold for the non-L\"obell polytope in Figure~\ref{fig:lobell-edge-child}. It is obtained from $L_6$ by splitting $U$ along an edge joining interior points of $U\cap X_0$ and $U\cap X_3$. The resulting pentagons are labeled $U_+$ and $U_-$ as shown.

	Label every edge $2$, take $(F_0,F_1,F_2,F_3)=(D,U_+,X_0,X_1)$, and order the remaining facets as follows:
	\[
		Y_0,Y_1,X_2,Y_2,X_3,Y_3,U_-,X_4,Y_4,X_5,Y_5.
	\]
	The initial facets satisfy~\eqref{hyp:initial}, and the following choices of triples verify~\eqref{hyp:extension}:
\[
				\begin{array}{c|ccc}
					T_j & P_j & Q_j & R_j \\
					\hline
					Y_0 & D & X_0 & X_1 \\
					Y_k\ (1\le k\le5) & D & X_k & Y_{k-1} \\
					X_k\ (k=2,3) & U_+ & X_{k-1} & Y_{k-1} \\
					U_- & U_+ & X_0 & X_3 \\
					X_k\ (k=4,5) & U_- & X_{k-1} & Y_{k-1}
				\end{array}
			\]
\begin{figure}[H]
\centering
\begin{tikzpicture}[x=0.0130cm,y=-0.0130cm,scale=0.7,transform shape,line width=0.45pt,line cap=round,line join=round,font=\normalsize,every node/.style={inner sep=0pt}]
				\path[use as bounding box] (20,20) rectangle (535,610);
				\foreach \name/\x/\y in {
					d0/279/35,d1/521/175,d2/521/454,d3/279/593,d4/36/454,d5/36/175,
					c0/148/238,c1/283/156,c2/414/231,c3/411/391,c4/274/473,c5/144/398,
					b0/235/219,b1/335/215,b2/385/307,b3/324/408,b4/222/414,b5/173/324,
					a0/270/260,a1/313/256,a2/330/296,a3/287/370,a4/245/373,a5/227/334,
					w0/264/307,w3/295/322
				} {
					\coordinate (\name) at (\x,\y);
				}
				\foreach \i in {0,...,5} {
					\pgfmathtruncatemacro{\j}{mod(\i+1,6)}
					\draw (d\i)--(d\j);
					\draw (d\i)--(c\j);
					\draw (c\i)--(b\i)--(c\j);
					\draw (b\i)--(a\i);
				}
				\draw (a5)--(w0)--(a0)--(a1)--(a2)--(w3)--(a3)--(a4)--(a5);
				\draw (w0)--(w3);
				\node at (470,86) {$D$};
				\node at (294.5,287.2) {$U_+$};
				\node at (263.4,341.9) {$U_-$};
				\foreach \i/\x/\y in {
					0/211.5/275.7,1/285.8/215.4,2/361.4/259.4,
					3/346.5/353.4,4/271.8/413.7,5/196.4/370.4
				} {
					\node at (\x,\y) {$X_{\i}$};
				}
				\foreach \i/\x/\y in {
					0/186.5/152.1,1/374.1/149.2,2/465.9/313.2,
					3/370.9/476.3,4/183.1/479.5,5/91.7/315.6
				} {
					\node at (\x,\y) {$Y_{\i}$};
				}
			\end{tikzpicture}
\caption{The $1$-skeleton of the unique edge-child of $L_6$, redrawn from \cite[Appendix~B, polyhedron~3]{Ino2015}. The unbounded face is $D$.}
		\label{fig:lobell-edge-child}
\end{figure}
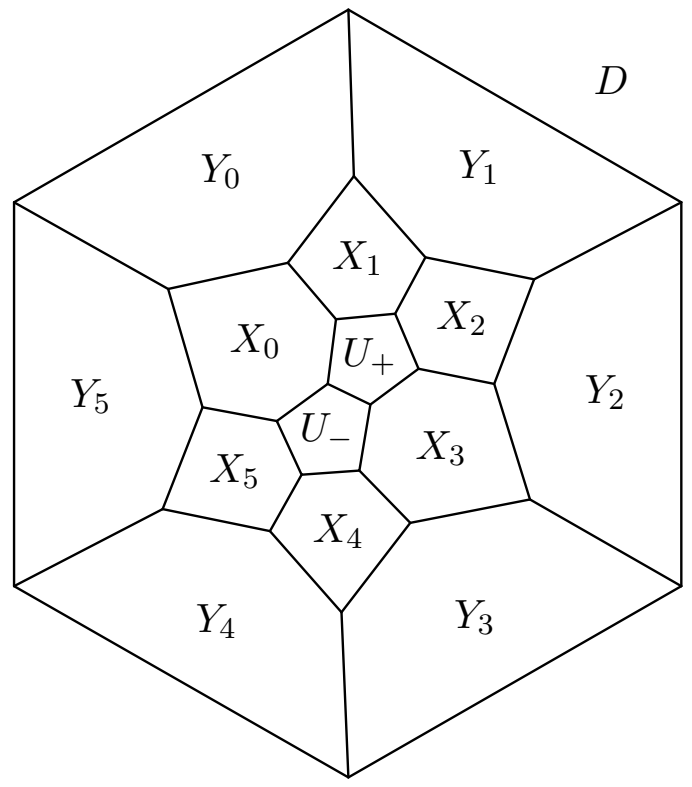

\bibliographystyle{amsplain}
\bibliography{references}

\end{document}